\documentclass[12pt,epsfig,tikz,multi]{amsart}
\usepackage{url}
\usepackage{verbatim}
\usepackage{latexsym}
\include{plaquettemacro}
\usepackage{hyperref}
\usepackage{amsrefs}
\usepackage{breqn}
\usepackage{float}
\restylefloat{table}
\usepackage{amsmath,amscd,amssymb}
\usepackage{graphicx}
\usepackage{caption}
\usepackage{subcaption}
\usepackage{enumitem}
\usepackage{xcolor}
\usepackage{graphicx,wrapfig,lipsum}
\usepackage{tikz}
\usepackage{diagbox}
\usepackage[utf8]{inputenc}
\usepackage[T1]{fontenc}
 2

\newtheorem{theorem}{Theorem}[section]

\newtheorem*{theorem*}{Theorem}
\newtheorem{definition}[theorem]{Definition}

\numberwithin{equation}{section}
\title[About cuspidal crosscaps and folded singularities]{Cuspidal crosscaps and folded singularities on a maxface and a minface} 
\author{Rivu Bardhan, Anu Dhochak, Pradip Kumar}
\address{Department of Mathematics, Shiv Nadar University, Dadri 201314, Uttar Pradesh, India.}
\email{rb212@snu.edu.in, ad404@snu.edu.in, pradip.kumar@snu.edu.in}
\date{}


\begin{document}
\maketitle
\begin{abstract}
 For a given zero mean curvature surface $X$ (in the Lorentz Minkowski space) having folded singularity,  we construct a family of maxface and minface, having increasing cuspidal crosscaps,  converging to $X$.  We include a general discussion of this. 
\end{abstract}

\section{Introduction}
Space-like (time-like) zero mean curvature surfaces in the Lorentz-Minkowski space $\mathbb E_1^3$ are locally area maximizing (minimizing respectively) surfaces. These are similar to the minimal surface in $\mathbb R^3$.  But in contrast to the minimal surface in $\mathbb R^3$, these surfaces appear with non isolated singularities \cite{shintaro1}, \cite{Estudillo1992}, \cite{fujimori2015zero}, \cite{KOBAYASHI1983}.   The generalized maximal (time-like minimal) immersion having only non-isolated singularities are called maxface (minface respectively).   

These non-isolated singularities are of various types.  As front or frontal, on these surfaces, we see singularities as cuspidal edge, cuspidal crosscaps, swallowtails, etc. There are other space-like (time-like) zero mean curvature surfaces that are  neither front nor  frontal but appear with other types of singularities, such as  folded singularities, cone-like, etc.

In \cite{Kim2006}, the authors gave an example of a family of maxface $X_n$ for each natural number $n$, that  has swallowtails in an increasing order.  Further, in  \cite{Fujimori2009}, the authors talked about Trinoids (a maxface) with swallowtails whose computer graphics look cone-like singularities.  Moreover, they  remarked that the computer graph of many swallowtails together looks cone-like.  In \cite{anupradip}, the authors justified this phenomenon by producing a sequence of maxface with an increasing number of swallowtails, converging to the maxface with cone-like singularity. 

We see similar phenomena (at least graphically) in cuspidal crosscaps. The graph of  cuspidal crosscaps looks like the folded singularity. Similar to swallowtails and cone-like discussed in \cite{anupradip}, we can ask a question:  given a maxface (minface) having folded singularity, does there exist a family of maxface (minface respectively) having an increasing number of cuspidal crosscaps and `converging' to the former one?

 In this article, in section 4, we give such a family and prove the following: 
\begin{theorem*}[Theorem \ref{thm:Last}]
Let X be a maxface (minface)  having folded singularity.  Then there is a sequence $X_n$ of maxfaces (sequence $Y_n$ of minfaces) having an increasing number of cuspidal crosscaps, and the sequence  `converges' (in the norm $\|.\|_\Omega$) to X.
\end{theorem*}

To construct such a family, we will find the condition on  the singular Bj\"{o}ring data to have a cuspidal crosscap.  For the maxface, such conditions are already derived in \cite{SaiPradip}.  In this article, in  section 3, we find similar conditions on the singular Bj\"{o}rling data for the minface, such that they have a cuspidal edge, swallowtails, and cuspidal crosscaps. 

Besides finding the sequence above, such conditions may be useful in  various situations.  For example, these are required initial conditions to find interpolating surfaces containing given two curves with the prescribed nature of singularities.

\section{Preliminaries}
The Lorentz Minkowski space is a vector space $\mathbb R^3$ with the metric $ds^2 := dx^2+dy^2-dz^2$; we denote it by $\mathbb E_1^3$.  Space-like zero mean curvature immersion is called maximal immersion. If we allow singularities, then we call it a generalized maximal immersion.  Moreover, if singularities are due to those points where the limiting tangent plane contains a light-like vector, then we call such a generalized maximal immersion as maxface \cite{UMEHARA2006}. In the following, we will discuss two methods to construct the generalized maximal immersion \cite{Estudillo1992,UMEHARA2006,RPR2016,Kim2007}. 
\subsection{Weierstrass-Enneper representation for the generalized maximal immersion} 
 Let $M^2$ be a Riemann surface.  Let $g$ be a meromorphic function  and $\omega$ a holomorphic 1-form on $M^2$ such that   $Re\int_\gamma \left(\frac{1}{2}(1+g^2), \frac{i}{2}(1-g^2), -g\right)\omega = 0$ for all loops $\gamma$ on $M^2$.  Then the map
$$f(z) := Re\int_{z_0}^{z}\left(\frac{1}{2}(1+g^2), \frac{i}{2}(1-g^2), -g\right)\omega$$
is well-defined and gives a generalized maximal immersion \cite{Estudillo1992,UMEHARA2006}. Moreover, any generalized maximal immersion can be written in this manner. 

\subsection{Singular Bj\"{o}rling representation for the generalized maximal immersion}\label{section:singular Bjorling Representation for maximal surface}
Let  $ \gamma : I \to \mathbb{E}_1^3$ be a real analytic null curve and  $L: I \to \mathbb{E}_1^3$ be a real analytic null vector field such that $\langle \gamma'(u) , L(u) \rangle = 0$ and  $\gamma' \not \equiv 0 \text{ or } L \not \equiv 0$. Let  $g(z)$ be the analytic extension of 
$$g(t):=\begin{cases}
\frac{\gamma_1^\prime(t)+i\gamma_2^\prime(t)}{\gamma_3^\prime(t)}\; \text{ if } \gamma_3^\prime(t)\not \equiv 0\\
\frac{L_1(t)+iL_2(t)}{L_3(t)} \text{ if } L_3(t) \not \equiv 0.
\end{cases}$$
defined on  $U \subset \mathbb{C}$, $I \subset U$. If $|g(z)| \not \equiv 1$,
then the map
\begin{equation}\label{eq:Singular Bjorling Representation for maximal surface}X_{max,\gamma, L} (u + iv) :=\frac{1}{2}(\gamma(u+iv)+ \gamma(u- iv))- \frac{i}{2}\int_{u-iv}^{u+iv} L(w) dw
\end{equation} is the generalized maximal surface with $u - v$ as a conformal parameter and $ X_{max,\gamma, L}(u,0) = \gamma(u),\; \frac{\partial }{\partial v}\lvert_{(u,0)}X_{max,\gamma, L} = L(u)$ \cite{RPR2016,Kim2007}.

Similar to the maximal immersion, we have time-like zero mean curvature immersion and minface. %are those time-like zmc immersions where at all singular point limiting tangent plane contains a light like vector \cite{shintaro1}. 
Below we will recall the Weierstrass Enneper representation and the  Bj\"{o}rling formula for a minface \cite{shintaro1, KYKSEYY2011}.  
\subsection{Weierstrass-Enneper representation for minface  \cite{shintaro1} }  Let $M^2$ be a smooth manifold. A smooth map $X: M^2 \to \mathbb{E}_1^3$ is  a minface if and only if  at each point of $M^2$, there exists a local coordinate system $(u, v)$ in a domain $U$, the functions $ g_1(u)$, $g_2(v)$, and the 1-forms $\omega_1 = \hat{\omega_1}du$, $\omega_2 = \hat{\omega_2}dv $ such that 
\begin{enumerate}
\item $g_1g_2 \not\equiv 1$ on an open dense subset of $U$
 \item $\hat{\omega_1} \neq 0$, $\hat{\omega_2} \neq 0$ at each point on $U$.

\end{enumerate} Moreover, $X$ can be decomposed into two regular null curves. 
\begin{equation}\label{eqn: WER for minface}
    X(u,v) = \frac{1}{2} \left(\int_{u_0}^{u} (2g_1,1-{g_1}^2,-1-{g_1}^2)\omega_1 + \int_{v_0}^{v} (-2g_2, 1-{g_2}^2, 1+{g_2}^2)\omega_2 \right) + X(u_0,v_0).
\end{equation} 
 The quadruple $(g_1, g_2, \omega_1, \omega_2)$ is called the real Weierstrass data. The set of singular points on $U$ of a  minface $X$ correspond to the set $\{(u,v) \in U:  g_1(u)g_2(v) = 1\}.$

\subsection{Singular Bj\"{o}rling representation for the generalized time-like minimal surface }\label{Section:Singular Bjorling Representation for minfaces}
Let $\gamma : (a,b) \to \mathbb{E}_1^3$   be a smooth null curve, and  $L:(a,b) \to \mathbb{E}_1^3$  be a smooth  null vector field such that $\langle \gamma'(t), L(t) \rangle \equiv 0$ and 
    $\gamma' \not \equiv 0 \text{ or } L \not \equiv 0 \text{ and } \gamma' \not \equiv \pm L$.
Then the map
\begin{equation} \label{eqn: singular represetation minimal}X_{min,\gamma, L}(u,v) := \frac{1}{2}( \gamma(u+v) + \gamma(u-v)) + \frac{1}{2} \int_{u-v}^{u+v} L(s)ds, \end{equation} defined 
on $\Diamond_{a,b}:= \{(u,v)\in\mathbb R^2: a<u+v<b,\text{ and } a< u-v<b\}$  is a generalized time-like minimal surface.
Furthermore, 
$ X_{min,\gamma, L}(u,0) = \gamma(u),\;\frac{\partial}{\partial v} X_{min,\gamma, L}(u,0) = L(u), \;\;\forall u\in (a,b).$

Both maxface and minface admit various singularities. In the table \eqref{tab1}, we mention  the criterion given in \cite{SaiPradip} on the singular Bj\"{o}rling data  for maxface having various singularities. But for the minface we do not have a criterion to check the type of singularity in terms of the Bj\"{o}rling data till now. 

\begin{table}[ht]
\small
\begin{tabular}{|c |c| c| c| c| c| c|c|c| }\hline
\diagbox[width=4cm]{{Nature of}\\ {singularities at $p$}}{{Function's value}\\ at $p$}&$\gamma^\prime$&$L$&$\gamma^{\prime\prime}$&$\gamma^{\prime\prime\prime}$&$L^\prime$&$L^{\prime\prime}$&$\gamma_1^{\prime}\gamma_2^{\prime\prime}-\gamma_1^{\prime\prime}\gamma_2^{\prime}$&$L_1L_2^{\prime}-L_1^{\prime}L_2$\\
\hline
Cuspidal-edge&$\neq0$&$\neq0$&--&--&--&--&$\neq 0$&$\neq 0$\\
\hline
Swallowtails&$=0$&$\neq0$&$\neq0$&--&--&--&--&$\neq 0$\\
\hline
Cuspidal-Crosscaps&$\neq0$&$=0$&--&--&$\neq 0$&--&$\neq 0$&--\\
\hline
\end{tabular}
\caption{Criterion on the Singular Bj\"{o}rling data for maxface}
\label{tab1}
\end{table}

For the minface, the criterion to check the type of singularity in terms of the Weierstrass data is given in Fact 4.1 of \cite{shintaro1}. We recall it here.  The minface with Weierstrass data $(g_1, \hat{\omega_1}, g_2, \hat{\omega_2})$ has 

\begin{enumerate}
    \item  Cuspidal edge at $p$ if and only if at $p$, 
    \begin{equation}\label{equ: WER conditions for CE}
        \frac{g_1'}{g_1^2\hat{\omega_1}}-\frac{g_2'}{g_2^2\hat{\omega_2}}\neq 0,\; \frac{g_1'}{g_1^2\hat{\omega_1}}+\frac{g_2'}{g_2^2\hat{\omega_2}}\neq 0.\end{equation}
       \item  Cuspidal crosscap at $p$ if and only if at $p$, 
\begin{equation} \label{equ: WER conditions for crosscap } \frac{g_1'}{g_1^2\hat{\omega_1}}-\frac{g_2'}{g_2^2\hat{\omega_2}}=0,\; \frac{g_1'}{g_1^2\hat{\omega_1}}+\frac{g_2'}{g_2^2\hat{\omega_2}}\neq 0,\;\text{ and } \left(\frac{g_1'}{g_1^2\hat{\omega_1}}\right)^\prime\frac{g_2^\prime}{g_2}+ \left(\frac{g_2'}{g_2^2\hat{\omega_2}}\right)^\prime \frac{g_1^\prime}{g_1}\neq 0.\end{equation}
\item  Swallowtails at $p$ if and only if at $p$, 
\begin{equation} \label{equ: WER conditions for swallowtails } \frac{g_1'}{g_1^2\hat{\omega_1}}-\frac{g_2'}{g_2^2\hat{\omega_2}}\neq0,\; \frac{g_1'}{g_1^2\hat{\omega_1}}+\frac{g_2'}{g_2^2\hat{\omega_2}}= 0,\;\text{ and } \left(\frac{g_1'}{g_1^2\hat{\omega_1}}\right)^\prime\frac{g_2^\prime}{g_2}+ \left(\frac{g_2'}{g_2^2\hat{\omega_2}}\right)^\prime \frac{g_1^\prime}{g_1}\neq 0.\end{equation}
\end{enumerate}

The above singularities of a maxface and minface appear when we consider these  as front or frontal. Now we discuss another type of singularity called folded singularity that appears on maxface and minface but not as frontal or front.   
\begin{definition}[Folded singularity]\label{defn:folded}\cite{Kim2007}.  The $X_{max, \gamma, L}$ and $X_{min, \gamma, L}$ have folded singularities along $\gamma$ if $L=0$ and $\gamma$ is non degenerate curve. 
\end{definition}

There are many other singularities, but we recalled only selected one, as per the required discussion of this article.  For more details on the topic, we refer \cite{PSV} and \cite{KYKSEYY2011}. 

\section{Conditions on the singular Bj\"{o}rling data for minface}
We start with $\{\gamma, L\}$ the singular Bj\"{o}rling data such that 
\begin{itemize}
    \item At those points where $\gamma^\prime\neq 0$, we have  $\gamma_2^\prime \neq \pm \gamma_3^\prime$,  and
    \item A those points where $L\neq 0$, we have  $L_2 \neq \pm L_3$.
\end{itemize}
Let $X_{min,\gamma, L}$ be the generalized time-like minimal surface as in \eqref{eqn: singular represetation minimal}. Below we will find the Weierstrass quadruple for $X_{min,\gamma L}$ in terms of $\gamma$ and $L$. 

\subsection{Relation between  the singular Bj\"{o}rling data and Weierstrass quadruple for the minface}
After changing the coordinate $(u,v)$ to $(t,s)$ by $s=u+v$ and $t=u-v$  and for some $x_0 \in (a,b) $, $X_{min,\gamma, L}$ as in \eqref{eqn: singular represetation minimal},  changes to
\begin{equation}\label{eqn: new variable time-like represnetation}
    \Tilde{X}(t,s) = \frac{1}{2} \left( \gamma(t) + \int_{t}^{x_0} L(\eta) d\eta \right) + \frac{1}{2} \left( \gamma(s) + \int_{x_0}^{s} L(\eta) d\eta \right).
\end{equation} 
We take  $ \alpha(t) := \frac{1}{2} \left(\gamma(t) + \int_{t}^{x_0} L(\eta) d\eta \right)$ and  $\beta(s) := \frac{1}{2} \left( \gamma(s) + \int_{x_0}^{s} L(\eta) d\eta \right)$.
 It turns out that both $\alpha$ and $\beta$ are null curves. 
$\alpha^\prime(t) = \frac{1}{2}(\gamma^\prime(t) - L(t))$, and $\beta^\prime(s) = \frac{1}{2}(\gamma^\prime(s) + L(s))$.\\
For those point $t$, where $\gamma^\prime \neq 0 $, we have  $L(t) = c(t)\gamma^\prime(t)$. This gives
\begin{align*}
  &&  \alpha_2^\prime(t) - \alpha_3^\prime(t) =\frac{1}{2}(\gamma_2^\prime(t) - \gamma_3^\prime(t))(1-c(t))\\
 &&\beta_2^\prime(s) + \beta_3^\prime(s)= \frac{1}{2}(\gamma_2^\prime(s) + \gamma_3^\prime(s))(1+ c(s)).
\end{align*}
Similarly at those points where $L \neq 0$, $\gamma^\prime(t) = d(t)L(t)$, we get
\begin{align*}
    &&\alpha_2^\prime(t) - \alpha_3^\prime(t) = \frac{1}{2}(L_2(t) - L_3(t))(-1+d(t))\\
&&\beta_2^\prime(s) + \beta_3^\prime(s) = \frac{1}{2}(L_2(s) + L_3(s))(1+ d(s)).
\end{align*}
The conditions on the singular Bj\"{o}rling data as mentioned in the starting  of this section and in the subsection \ref{Section:Singular Bjorling Representation for minfaces}, imply that for all points $\alpha_2^\prime-\alpha_3^\prime\neq 0$, $\beta_2^\prime+\beta_3^\prime\neq 0$.

Let $\left(g_1,g_2,\hat{\omega_1},\hat{\omega_2} \right)$ be the Weierstrass quadruple for $\Tilde{X}$.  $\Tilde{X}$ has the representation as in  \eqref{eqn: WER for minface}.   By  comparing the two representations \eqref{eqn: WER for minface} and \eqref{eqn: new variable time-like represnetation}, we get components of $\alpha^\prime$ and $\beta^\prime$ as the following: 
\begin{equation}\label{eqn: alpha and g w relations}
    \alpha_1'(u) = 2g_1(u)\hat{\omega_1}(u),\;\;\alpha_2'(u)= (1-{g_1}^2(u))\hat{\omega_1}(u),\;\; \alpha_3'(u)= (-1-{g_1}^2(u))\hat{\omega_1}(u)
\end{equation}
\begin{equation}
     \beta_1'(v)= -2g_2(v)\hat{\omega_2}(v),\;\; \beta_2'(v)= (1-{g_2}^2(v))\hat{\omega_2}(v),\;\; \beta_3'(v) =(1+{g_2}^2(v))\hat{\omega_2}(v).
\end{equation}

From  \eqref{eqn: alpha and g w relations} and since $\hat{\omega_1} \neq 0$, we get $g_1 = \frac{\alpha_1^\prime}{2\hat{\omega_1}}$. Also  $\alpha_2^\prime = (1-g_1^2)\hat{\omega_1}$.  Solving this, we get $\hat{\omega_1} = \frac{\alpha_2^\prime \pm \alpha_3^\prime}{2}$ and then  $g_1 = \frac{\alpha_1^\prime}{\alpha_2^\prime \pm \alpha_3^\prime}$.

Now   for the case when $g_1 = \frac{\alpha_1^\prime}{\alpha_2^\prime - \alpha_3^\prime}$ and $\hat{\omega_1} = \frac{\alpha_2^\prime - \alpha_3^\prime}{2}$,  we get
\[-(1+g_1^2)\hat{\omega_1} = \frac{-1}{2} \left[\frac{\alpha_2^{\prime 2} + \alpha_3^{\prime 2} - 2 \alpha_2^\prime\alpha_3^\prime + \alpha_1^{\prime 2}}{\alpha_2^\prime - \alpha_3^\prime}\right] = \alpha_3^\prime.\]

For other choices of $\hat{\omega_1}, g_1$, we see that $-(1+g_1^2)\hat{\omega_1}\neq\alpha_3^\prime$.   Therefore we must have $$ g_1 = \frac{\alpha_1^\prime}{\alpha_2^\prime - \alpha_3^\prime}\; \text{and}\; \hat{\omega_1} = \frac{\alpha_2^\prime - \alpha_3^\prime}{2}. $$

Similarly we have $$g_2(v) = \frac{-\beta_1'}{\beta_2'+ \beta_3'}\; \text{and} \; \hat{\omega_2}(v) = \frac{\beta_2' + \beta_3'}{2}.$$

Now to get the Weierstrass quadruple in terms of the singular Bj\"{o}rling data, we replace the values of $\alpha$ and $\beta$ in terms of $\gamma$ and $L$.  We will do this case by case.

 If at $u$,  $\gamma^\prime \neq 0$, then $L(u)= c(u)\gamma^\prime(u)$ for the real function $c = \frac{L_3}{\gamma_3'}$ and we get, 
 \begin{equation}\label{gamma not zero}
  \hat{\omega_1} = \frac{\gamma_2' - \gamma_3'}{4}(1 - c) ,\;\; g_1 = -\frac{\gamma_2' + \gamma_3'}{\gamma_1'},\;\;   
    \hat{\omega_2}= \frac{\gamma_2'+ \gamma_3'}{4}(1+c ) ,\;\;  g_2 = \frac{\gamma_2' - \gamma_3'}{\gamma_1'}.   
   \end{equation} 
We recall here the fact that at points where $\gamma^\prime\neq 0$, we have $\gamma_3^\prime \neq 0$, and as $\gamma_2^\prime\neq \gamma_3^\prime$, and $\gamma$ a null curve, this gives $\gamma_1^\prime\neq 0$. 

Similarly when  $L\neq 0$, with $d = \frac{\gamma_3^\prime}{L_3}$, we have 
  \begin{equation}\label{ L not zero}
   \hat{\omega_1}= \frac{L_2-L_3}{4}(-1+d),\;\; g_1 = -\frac{L_2 + L_3}{L_1},\;\;
    \hat{\omega_2}= \frac{L_2 + L_3}{4}(1+d),\;\; g_2 = \frac{L_2 - L_3}{L_1}.
\end{equation}  

Thus, we found the Weierstrass quadruple of $X_{min, \gamma, L}$ in terms of $\{\gamma, L\}$. In the next subsection we will convert conditions in terms of the singular Bj\"{o}rling data so that a time-like minimal surface has various singularities. 
  \subsection{Conditions on the singular Bj\"{o}rling data of minface for various singularity} 
  Let $\lbrace\gamma\text{, }L\rbrace$ be the singular Bj\"{o}rling data for the generalized time-like minimal surface $X_{min,\gamma,L}$,  as defined in the section \ref{Section:Singular Bjorling Representation for minfaces}. We denote 
 $$ \delta_1 := \frac{g_1'}{g_1^2\hat{\omega_1}} \text{ , }  \delta_2 := \frac{g_2'}{g_2^2 \hat{\omega_2}} \text{ , } D(\gamma_{23}^\prime , \gamma_{23}^{\prime\prime}) := \gamma_2^\prime \gamma_3 ^{\prime\prime} - \gamma_3 ^\prime \gamma_2^{\prime\prime} \text{ , } D(L_{23} L_{23}^ \prime):= L_2 L_3^ \prime - L_3 L_2 ^\prime.$$
 
We again start with the case when $\gamma^\prime \neq 0$ and $L= c\gamma^\prime$. Using \eqref{gamma not zero} we have 
\begin{align*}\delta_1 &= \frac{4g_1^\prime}{g_1^2(1-c)(\gamma_2^\prime - \gamma_3^\prime)}\\
        &= \frac{4[(\gamma_2^{\prime\prime} + \gamma_3^{\prime\prime})\gamma_1^\prime - \gamma_1^{\prime\prime}(\gamma_2^\prime + \gamma_3^\prime)]}{(\gamma_1^{\prime 2})(1 - c)(\gamma_2^\prime + \gamma_3^\prime)}\\
        &=\frac{4}{\gamma_1^{\prime 3}(1-c)}\left[-\gamma_1^{\prime\prime}\gamma_1^\prime - \gamma_2^{\prime\prime}\gamma_2^\prime + \gamma_3^{\prime\prime}\gamma_3^\prime + \gamma_2^\prime\gamma_3^{\prime\prime} - \gamma_3^\prime\gamma_3^{\prime\prime}\right]\\
        &=\frac{-4}{\gamma_1^{\prime 3}(1-c)}\left[\langle \gamma^\prime, \gamma^{\prime\prime}\rangle + D(\gamma_{23}^\prime, \gamma_{23}^{\prime\prime}) \right].
    \end{align*}
    Since $\gamma$ is a singular curve, $\langle \gamma^\prime, \gamma^{\prime\prime} \rangle = 0$. This gives $ \delta_1=-4\frac{D(\gamma_{23}^\prime,\gamma_{23}^{\prime\prime})}{\gamma_1^{\prime 3}(1-c)}$.
 
    Moreover at those points where $c(u)=0$, we have  
    $$
   \delta_1^\prime = -4\frac{\gamma_1^\prime D(\gamma_{23}^\prime,\gamma_{23}^{\prime\prime\prime})- D(\gamma_{23}^\prime,\gamma_{23}^{\prime\prime})(3\gamma_1^{\prime\prime} - \gamma_1^\prime c')}{\gamma_1^{\prime 4}}. 
$$
    A similar calculation gives
    $\delta_2 = -4\frac{D(\gamma_{23}^\prime,\gamma_{23}^{\prime\prime})}{\gamma_1^{\prime 3} (1+c)}$, and at those points where $c(u)=0$,$$ \delta_2^\prime= -4\frac{\gamma_1^\prime D(\gamma_{23}^\prime,\gamma_{23}^{\prime\prime\prime})- D(\gamma_{23}^\prime,\gamma_{23}^{\prime\prime})(3\gamma_1^{\prime\prime} + \gamma_1^\prime c')}{\gamma_1^{\prime 4}}.$$
    As $c(u)= \frac{L_3(u)}{\gamma_3^\prime(u)}$, we get
    \begin{align*}
     \delta_1 + \delta_2 &= \frac{g_1'}{g_1^2 \hat{\omega_1}} + \frac{g_2'}{g_2^2 \hat{\omega_2}}\\ 
      &=\frac{-8}{\gamma_1^{\prime 3} (1 - c^2)} \left[ \langle \gamma', \gamma ^{\prime\prime} \rangle c  + D(\gamma_{23}^\prime ,\gamma_{23}^{\prime\prime}) \right]\\
      &= \frac{-8\gamma_3^\prime}{\gamma_1^{\prime 3}(\gamma_3^{\prime 2} - L_3^{\prime 2})} \left[ \langle\gamma', \gamma ^{\prime\prime} \rangle L_3 + \gamma_3^\prime D(\gamma_{23}^\prime, \gamma_{23}^{\prime\prime} )\right]\\
      &= \frac{-8\gamma_3^{\prime 2} D(\gamma_{23}^\prime, \gamma_{23}^{\prime\prime} )}{\gamma_1^{\prime 3}(\gamma_3^{\prime 2} - L_3^{\prime 2})}, \text{ and }\\
      \delta_1 - \delta_2 &= \frac{g_1^\prime}{g_1^2\hat{\omega_1}} - \frac{g_2^\prime}{g_2^2\hat{\omega_2}} = \frac{-8 \gamma_3^\prime L_3 D(\gamma_{23}^{\prime} , \gamma_{23}^{\prime\prime})}{\gamma_1^{\prime 3} (\gamma_3^{\prime 2} - L_3^{\prime 2})}.
  \end{align*}

Below we find  the conditions for a generalized time-like minimal surface having singularities like cuspidal edge, cuspidal crosscaps and swallowtails in term of singular Bj\"{o}rling data $\lbrace\gamma, L\rbrace$.
 \subsubsection{Cuspidal edge}Using  \eqref{equ: WER conditions for CE},  the time-like minimal surface $X_{min,\gamma, L }$ for  the data $\{\gamma, L\}$ has cuspidal edge at $u$  if and only if $\gamma^\prime \neq 0$, $L\neq 0$, $D( \gamma_{23}^\prime, \gamma_{23}^{\prime\prime}) \neq 0$.
 \subsubsection{Cuspidal crosscaps} At those point where $\gamma^\prime\neq 0$,  
 $ \delta_1 - \delta_2 =0, \text{ and } \delta_1 + \delta_2 \neq 0$ give 
 \begin{equation}\label{eqn:condition for cuspidal crosscaps1}
     \gamma_3^\prime \neq 0,\;\; L_3 = 0,\;\; D( \gamma_{23}^\prime, \gamma_{23}^{\prime\prime}) \neq 0.  
 \end{equation}
     
     Moreover at those point where \eqref{eqn:condition for cuspidal crosscaps1} satisfies, we have
     $$\nonumber\delta_1'\frac{g_2'}{g_2} + \delta_2^\prime \frac{g_1'}{g_1} = \frac{g_1'}{g_1}(\delta_2^\prime - \delta_1^\prime)
         = -\frac{4}{\gamma_1^{\prime 6}}D(\gamma_{23}^\prime,\gamma_{23}^{\prime\prime})\left[ D(\gamma_{23}^\prime,\gamma_{23}^{\prime\prime})(2c' \gamma_1^\prime) \right]$$
         and
         $c^\prime = \frac{L_3^\prime \gamma_3^\prime - L_3 \gamma_3^{\prime\prime}}{\gamma_3^{\prime 2}}.$
         Therefore $$\delta_1'\frac{g_2'}{g_2} + \delta_2^\prime \frac{g_1'}{g_1} \neq 0 \iff c' \neq 0\iff L^\prime \neq 0.$$
   Using  \eqref{equ: WER conditions for crosscap },  at $u$, where $\gamma^\prime\neq 0$,  the time-like minimal surface $X_{min,\gamma, L }$ for  the data $\{\gamma, L\}$ has cuspidal crosscap if and only if 
 $\gamma^\prime\neq 0,\;\; L=0, \;\; L^\prime\neq 0, \text{ and } D(\gamma_{23}^\prime, \gamma_{23}^{\prime\prime}) \neq 0.$
 \subsubsection{Swallowtails} Similar to the cuspidal crossscaps,   at $u$, where $\gamma^\prime\neq 0$,  the time-like minimal surface $X_{min,\gamma, L }$ for  the data $\{\gamma, L\}$ has swallowtails if and only if 
 $\gamma^\prime(u)= 0,\;\; L(u)\neq 0, \;\; \gamma^{\prime \prime}(u)\neq 0 \text{ and } D(L_{23}, L_{23}^{\prime}) \neq 0.$

  For a non-degenerate curve (the curve is  non-degenerate if and only if for all $u\in I$, $D(\gamma_{23}^\prime,\gamma_{23}^{\prime\prime})\neq 0$ or $D(\gamma_{12}^{\prime},\gamma_{12}^{\prime\prime})\neq 0$), these conditions turn out to be similar to the one given in the table \eqref{tab1} for the maxface. 
 
%From the above conditions and the table \eqref{tab1},  below we construct the family of minface and maxface having cuspidal crosscap singularities. Though we have already mentioned the purpose of taking the cuspidal crosscaps, we write it here.  If we choose the sequence of data $\lbrace\gamma, L_n\rbrace$ with increasing order of cuspidal edge, we see that it `seems' to converge to the folded one, and folded singularity is the one that gives a chance to patch space-like and time-like surfaces. So mix-type surfaces can be viewed as a limit of maxface and minface having cuspidal crosscaps.   Below we will make this more precise. 

\section{Family of maxfaces and minfaces}\label{section:family of maxfaces and minfaces joined along cuspidal edge}
We start this section with an example: Helicoid. This is an example of how the time-like and space-like Helicoid patches smoothly along the folded singularity. For details, we refer to \cite{KRSUY2005}.

\subsection{Mix-type Helicoid and a sequence} Let  $\gamma: I\to \mathbb{E}_1^3;\;\;\gamma(t):=(\cos{t},\sin{t},t)$.     With $L=0$, $\lbrace\gamma,L\rbrace$ turns out to be the singular Bj\"{o}rling data for both generalized time-like minimal surface and generalized space-like maximal surface.

For this data,  we have  the following generalized space-like maximal surface
\begin{equation}\label{sph}
  X_{max,\gamma,L}(u, v)=(\cos{u}\cosh{v},\sin{u}\cosh{v},u)  
\end{equation}
on $\Omega$ and  $\gamma$ is  its singular curve corresponding to the points $v=0$.
 Also, for the same data, we have a generalized time-like minimal surface as following 
\begin{equation}\label{tmh}
    X_{min,\gamma,L}(u,v)=(\cos u\cos v,\sin u \cos v,u).
\end{equation}
The generalized immersions as in   \eqref{sph}, \eqref{tmh} are  parametrizations of the helicoid of type space-like and time-like, respectively.  We see that both these surfaces contain $\gamma$ and have folded singularity on $I$.

Without loss of generality we assume that $I= (0,2\pi)$ and let $r>0$ be such that $\overline{\Omega}, \overline{\Diamond}\subset B(0, r)$.   With this $r$ and $n\geq 1$, we take the singular Bj\"{o}rling data as 

\begin{equation}\label{eqn:seqbjorlingdata}\left\{ \gamma,\gamma^\prime\frac{1}{n(r+1)^n}\prod_{i=1}^{n}(w-\frac{1}{i})\right\}.
\end{equation}  For each $n\in\mathbb{N}$, let $X_{max,\gamma, L,n}$, be the generalized maximal immersion as in  \eqref{eq:Singular Bjorling Representation for maximal surface} and $X_{min,\gamma, L,n}$  be the generalized time-like minimal immersion as in \eqref{eqn: singular represetation minimal}.  This defines a pair $\{X_{max,\gamma, L,n}, X_{min, \gamma, L, n}\}$ of generalized space-like maximal surfaces and the generalized time-like minimal surface containing $\gamma$ as a singular curve at $v=0$.

Using the condition we found in section 3, it is easy to verify that both the generalized maximal immersions and generalized minimal immersions $\{X_{max,\gamma, L, n}, X_{min,\gamma, L, n}\}$  has $n$ cuspidal crosscaps at $t= \frac{1}{k}$. Moreover corresponding cuspidal crosscaps are  $\lbrace(\cos{\frac{1}{k}},\sin{\frac{1}{k}},\frac{1}{k})\rbrace_{k=1}^{n}\subset\gamma(I)$.

Helicoid is an example of a mix-type zero mean curvature surface, where the casual character changes along the folded singularity (see  \cite{KYKSEYY2011}). We have generated a sequence of immersions  $\lbrace X_{max,\gamma, L, n}, X_{min,\gamma, L, n}\rbrace_{n\in\mathbb{N}}$ with the increasing number of cuspidal crosscaps. For each $n$, $X_{max,\gamma, L, n}$ and $X_{min,\gamma, L, n}$  are joined along $\gamma$. Moreover the sequence $\lbrace X_{max,\gamma,L,n}\rbrace$ `converges' towards the mix-type helicoid (that is smoothly patched along a folded singularity) as $n\rightarrow\infty$.

In the next subsection, we will start with a generalized mixed-type surface that is smoothly joined along a folded singularity (for example: helicoid).  Next, we will find a sequence of mixed-type surfaces joined along a null curve with increasing cuspidal crosscaps and this sequence `converges' to the surface we started with.  

\subsection{Family for a general mix-type zero mean curvature  surfaces}

Authors in \cite{anupradip,SaiPradip} have defined and discussed the convergence of maxfaces in a particular norm. We would recall the norm here. 
  Let $\Omega\subset \mathbb C$ be a bounded simply connected domain, $\overline{\Omega}$ be its closure. Let $X\in C(\overline{\Omega},\mathbb R^3)$, the space of continuous maps. For each $z \in\overline{\Omega}$, we define 
$$
\|X(z)\|:={\rm max}\{\lvert X_1(z)\rvert,\lvert X_2(z)\rvert,\lvert X_3(z)\rvert\}\,{ \rm and  } \, \|X\|_{\Omega}:=\sup_{z\in{\overline \Omega}}\|X(z)\|.
$$

We start with the singular Bj\"{o}rling data $\{\gamma, L=0\}$, defined on a closed and bounded interval $I$, without loss of generality we will assume $I= [0,1]$. Moreover,  $X_{max,\gamma, 0}$ and $X_{min,\gamma, 0}$ be the generalized maximal immersion  and generalized time-like minimal immersion, respectively, as in  \eqref{eq:Singular Bjorling Representation for maximal surface} and \eqref{eqn: singular represetation minimal} defined on some bounded open subsets $\Omega_1$ and $\Omega_2$ (of $\mathbb{R}^2$) that contains $I$. We also assume both have folded singularity along $\gamma$. Since $L=0$, $X_{max,\gamma,0}$ and $X_{min,\gamma, 0}$ have folded singularity at $v=0$. Therefore $\gamma$ is non-degenerate (see definition \ref{defn:folded} and the curve is  non-degenerate if and only if for all $u\in I$, $D(\gamma_{23}^\prime,\gamma_{23}^{\prime\prime})\neq 0$ or $D(\gamma_{12}^{\prime},\gamma_{12}^{\prime\prime})\neq 0$). 
%\begin{comment}
%\footnote{It is known that a topologically open set of %$\mathbb{C}$ is homeomorphic to an open set of $\mathbb{R}^2$ and vice versa through the map $\phi(x+iy):=(x,y)$. Therefore, the metric induced from the restriction of the metric of $\mathbb{R}^2$ on $\Omega_1$ can be induced on $\phi^{-1}(\Omega_1)$. With a little abuse of notation, we choose to call $\phi^{-1}(\Omega)$ with the metric induced from $\Omega_1$ as $\Omega_1$ itself. Therefore here we can consider $\Omega_1$ as a subset of $\mathbb{C}$ as per our need and define a maxface with respective singular Bj\"{o}rling data as provided because while calculating the norm on $C(\overline{\Omega},\mathbb{R}^3)$, we are not concerned about the differential structure on $\Omega_1$ but about the metric on $\Omega_1$ as a subset of $\mathbb{R}^2$.}
%\end{comment}

As both $\Omega_1$ and $\Omega_2$  are bounded  and $I= [0,1]\subset \Omega_1\cap \Omega_2$, there exists some $r>0$ and $\Omega$, $\Diamond$ such that $I\subset \Omega\subset \overline{\Omega_1}\subset B(0,r)$, and $I\subset \Diamond\subset \overline{\Omega_2}\subset B(0, r)$.   With this $r$ and $n\geq 1$, we take  the singular Bj\"{o}rling data as in  \eqref{eqn:seqbjorlingdata}.  For this data, we have corresponding immersions as follows. 
\begin{equation}\label{equation:cuspidal crosscap maximal surface}
    X_{max,\gamma, L, n}(u, v)=\frac{1}{2}(\gamma(u+iv)+\gamma(u-iv))-\frac{i}{2n(r+1)^n}\int_{u-iv}^{u+iv}\gamma^\prime(s)\prod_{k=1}^{n}(s-\frac{1}{k})ds
\end{equation}
on $\Omega_1$.
\begin{equation}\label{equation: cuspidal crosscap minimal surface}
    X_{min,\gamma, L, n}(u,v)=\frac{1}{2}(\gamma(u+v)+\gamma(u-v))+\frac{1}{2n(r+1)^n}\int_{u-v}^{u+v}\gamma^\prime(s)\prod_{k=1}^{n}(s-\frac{1}{k})ds
\end{equation}
on $\Omega_2$.

Both the generalized immersions have  cuspidal-crosscap singularities at $\{\frac{1}{k}: k\in \{1,2,...,n\}\}$ for each $n$.   Moreover following the similar calculations as in \cite{anupradip,SaiPradip}, we see that 
\begin{align*}
    &\lVert X_{max,\gamma, L,n} - X_{max,\gamma, 0} \rVert_{\Omega_1} \leq M (\lVert \gamma_n^{k\prime}- \gamma^{k\prime} \rVert_{\Omega_1} + \lVert L_n^k - L^k \rVert_{\Omega_1} ) \leq M \lVert L_n \rVert_{\Omega_1} \to zero.\\
    &\text{Similarly, we have } \lVert X_{min, \gamma, L, n} - X_{min, \gamma, 0} \rVert_{\Omega_2} \leq M' \lVert L_n \rVert_{\Omega_2} \to zero.
\end{align*}  

Thus we summarize all our discussions in this section in the following theorem.
   \begin{theorem}\label{thm:Last}
   Let $\{X_{max}, X_{min}\}$ be a maxface and a minface having a folded singularity along $\gamma$.  Then there exists two sequences $X_{max, \gamma, L, n}$ and $X_{min, \gamma, L, n}$ of maxfaces and minfaces respectively such that 
   \begin{enumerate}
       \item On the trace of $\gamma$ both immersions have an increasing number of cuspidal crosscaps at the same points. 
       \item $X_{max,\gamma, L,n}$ converges to $X_{max}$ in $\lVert . \rVert_{\Omega_1}$ and $X_{min,\gamma, L,n}$ converges to $X_{min}$ in $\lVert . \rVert_{\Omega_2}$.
   \end{enumerate}
   \end{theorem}

\medskip

\bibliography{ref}

\end{document}